\documentclass[11pt]{article}

\usepackage{amsmath,amsfonts,amsthm,amssymb,graphicx,tikz,tikz-cd}
\usepackage[all,2cell,ps]{xy}

\theoremstyle{plain}
\newtheorem{thm}{Theorem}[section]

\newtheorem{lem}[thm]{Lemma}
\newtheorem{prop}[thm]{Proposition}

\newtheorem{cor}[thm]{Corollary}

\newtheorem{qtn}[thm]{Question}

\theoremstyle{definition}

\newtheorem{rem}[thm]{Remark}

\theoremstyle{remark}

\newcommand{\bbB}{\mathbb{B}}
\newcommand{\bbC}{\mathbb{C}}

\newcommand{\bbP}{\mathbb{P}}
\newcommand{\bbQ}{\mathbb{Q}}

\newcommand{\bbZ}{\mathbb{Z}}

\newcommand{\calD}{\mathcal{D}}

\newcommand{\calO}{\mathcal{O}}

\newcommand{\calT}{\mathcal{T}}

\newcommand{\calW}{\mathcal{W}}

\newcommand{\al}{\alpha}
\newcommand{\gam}{\gamma}
\newcommand{\Gam}{\Gamma}

\newcommand{\de}{\delta}
\newcommand{\Del}{\Delta}

\newcommand{\Lam}{\Lambda}
\newcommand{\sig}{\sigma}
\newcommand{\Sig}{\Sigma}

\DeclareMathOperator{\U}{U}

\DeclareMathOperator{\SL}{SL}

\DeclareMathOperator{\PU}{PU}

\DeclareMathOperator{\Id}{Id}

\DeclareMathOperator{\Aut}{Aut}

\newcommand{\ssm}{\smallsetminus}

\newcommand{\wh}{\widehat}
\newcommand{\wt}{\widetilde}

\newenvironment{pf}{\begin{proof}}{\end{proof}}

\newenvironment{enum}{\begin{enumerate}}{\end{enumerate}}

\makeatletter
\let\@@pmod\pmod
\DeclareRobustCommand{\pmod}{\@ifstar\@pmods\@@pmod}
\def\@pmods#1{\mkern4mu({\operator@font mod}\mkern 6mu#1)}
\makeatother

\def\O{{\mathcal{O}}}
\def\L{{\mathcal{L}}}
\def\spec{\textrm{Spec}}

\usetikzlibrary{arrows}

\title{Rigid surfaces arbitrarily close to the Bogomolov--Miyaoka--Yau line}
\author{Matthew Stover\footnote{This material is based upon work supported by Grant Number 523197 from the Simons Foundation/SFARI and Grant Number DMS-1906088 from the National Science Foundation.} \\ \small{Temple University} \\ \small{\textsf{mstover@temple.edu}}
\and
Giancarlo Urz\'ua \footnote{The second author  was supported by the FONDECYT regular grant 1190066.} \\ \small{Pontificia Universidad Cat\'olica de Chile} \\ \small{\textsf{urzua@mat.uc.cl}}}
\date{\today}

\begin{document}

\maketitle

\begin{abstract}
We prove the existence of rigid compact complex surfaces of general type whose Chern slopes are arbitrarily close to the Bogomolov--Miyaoka--Yau bound of $3$. In addition, each of these surfaces has first Betti number equal to $4$.
\end{abstract}

\section{Introduction}\label{sec:Intro}

The list of complex surfaces of general type that are rigid is very short. From the point of view of geography, we know that their Chern slope $\frac{c_1^2}{c_2}$ can only take values in the interval $[\frac{5}{7},3]$  (see \cite[Thm.\ 2.3]{BauerCatanese}). A catalogue was recently collected in \cite{BauerCatanese}, where new families of rigid surfaces were also given. Together with the new interesting examples of rigid but not infinitesimally rigid surfaces in \cite{BauerPignatelli}, we can give the list according to their Chern slope as follows (for a more complete account, see the references in \cite{BauerCatanese}):

\begin{enum}
\item \textit{Ball quotients}: The universal cover is the two-dimensional ball, and they all satisfy $\frac{c_1^2}{c_2} =3$.

\item \textit{Irreducible bi-disk quotients}: The universal cover is the product of two upper half-planes, and they all satisfy $\frac{c_1^2}{c_2} =2$.

\item \textit{Beauville surfaces}: These are rigid unramified quotients of products of curves. They satisfy $\frac{c_1^2}{c_2} =2$.

\item \textit{Mostow--Siu surfaces}: These surfaces give an infinite sequence of slopes in $]2,2.9416]$ with accumulation point $\frac{8}{3}$. See Deraux \cite{Deraux} for more examples beyond the original examples of Mostow--Siu \cite{MostowSiu}.

\item \textit{Certain Catanese--Rollenske Kodaira fibrations}: Finitely many examples exist, and they satisfy $\frac{c_1^2}{c_2} = \frac{8}{3}$.

\item \textit{Certain Hirzebruch--Kummer coverings} \cite{Zuo}, \cite{BauerCatanese}: These are certain Galois coverings of the projective plane branched along rigid line arrangements, namely the complete quadrilateral, Ceva arrangement, and extended Ceva arrangement. These surfaces give an infinite discrete sequence of $\frac{c_1^2}{c_2}$ with infinitely many limit points in $[\frac{5}{2},\frac{8}{3}[$.

\item \textit{Bauer--Pignatelli surfaces} \cite{BauerPignatelli}: These are minimal resolutions of certain nodal product-quotient surfaces, which turn out to be rigid but not infinitesimally rigid. They give an infinite sequence $S_n$ with $\frac{c_1^2(S_n)}{c_2(S_n)} <2$ and $\lim_{n \to \infty}  \frac{c_1^2(S_n)}{c_2(S_n)}=2$.
\end{enum}

In particular, the above classes contain no rigid surfaces with Chern slope arbitrarily close to the Bogomolov--Miyaoka--Yau line. In this paper, we prove existence of such surfaces.

\begin{thm}\label{thm:Main}
For each odd integer $n \ge 3$, there exists a smooth minimal complex projective surface $V_n$ of general type such that:
\begin{enum}

\item $V_n$ has irregularity two;

\item $\frac{c_1^2(V_n)}{c_2(V_n)} = 3-\frac{4}{n^2}$;


\item $V_n$ has ample canonical divisor;

\item $V_n$ is rigid.

\end{enum}
\end{thm}

This is a theorem of Zuo in the case $n = 3$ \cite{Zuo}. In \cite[\S 2]{Zuo}, Zuo computes the irregularity of $V_3$ via completely different methods. In particular, our approach restricted to the case $n=3$ confirms Zuo's results.

We now briefly describe the construction of $V_n$. For $\zeta = e^{\pi i / 3}$, let $T$ be the elliptic curve $\bbC / \bbZ[\zeta]$ and $A$ be the abelian surface $T \times T$ with coordinates $(z, w)$. For $\al \in \{0, 1, \zeta\}$, if $T_\al$ denotes the curve $w = \al z$ on $A$, $T_\al$ is the curve $z = 0$ for $\al = \infty$, and $U_n$ is the group of $n$-division points on $A$, then $D_\al^{(n)} = U_n(T_\al)$ is a divisor with support $n^2$ disjoint smooth irreducible curves. When $n$ is odd, we take $V_n$ to be the smooth resolution of a certain $n$-fold cyclic cover of $A$ branched over the divisor $\sum D_\al^{(n)}$. See \S\ref{sec:Cyclic} for the precise definition of the branched cover.

The proof of Theorem \ref{thm:Main} exploits the fact that there is a torsion-free lattice $\Del_n < \PU(2,1)$ so that $\bbB^2 / \Del_n$ admits a smooth toroidal compactification biholomorphic to $V_n$, where $\bbB^2$ denotes the unit ball in $\bbC^n$ with its Bergman metric. Our construction is such that each $\Del_n$ is a normal subgroup of a fixed lattice $\Gam_1 < \PU(2,1)$ for which $\bbB^2 / \Gam_1$ has smooth toroidal compactification biholomorphic to the blowup of the abelian surface $A$ at $(0,0)$. We then use the structure of the group $\Gam_1$, in part aided by the computer algebra program Magma \cite{Magma}, to prove that $V_n$ has the listed properties. The definition of $\Gam_1$ and an account of some of its key properties is contained in \S\ref{sec:Hirz}. We define $V_n$ and prove all of Theorem \ref{thm:Main} except rigidity in \S\ref{sec:Cyclic}. Specifically, we prove (2) in Lemma \ref{lem:AmpleK}, (3) in Corollary \ref{cor:Irregularity},
and (4) along with the fact that $V_n$ is minimal of general type in Proposition \ref{prop:ChernCalc}.

To obtain that $V_n$ is rigid, the key step is to prove that $V_n$ has irregularity two. In fact, we show that $H^1(V_n,\Omega_{V_n}^1 \otimes \Omega_{V_n}^2)=0$ if $H^1(V_n,\O_{V_n})$ has dimension $2$. To prove rigidity from there, we consider an algebraic model of the surface $V_n$ as in \cite{EsnaultViehweg}, and use as main tools the sheaf of logarithmic differentials with poles along $\sum D_\al^{(n)}$ and a rigidity theorem of Fujiki for open ball quotients \cite{Fujiki}. The proof of rigidity is contained in \S\ref{sec:Rigidity}.

The surfaces $V_n$ and lattices $\Del_n$ also make sense when $n$ is even, but then $\Del_n$ is not normal in $\Gam_1$ and we exploit the fact that $\Del_n$ is normal in $\Gam_1$ when $n$ is odd in the proof that $V_n$ has irregularity two. Clearly the sequence $\{V_n\}$ for $n$ odd suffices to prove Theorem \ref{thm:Main}, hence we do not consider the case where $n$ is even. We suspect that $V_n$ is rigid for all $n$.

The construction of $V_n$ is very much inspired by Hirzebruch's construction of surfaces with $\frac{c_1^2}{c_2} \to 3$ using covers of the abelian surface $A$ \cite{Hirzebruch}. In fact, Hirzebruch notes in his paper that his examples are related to certain ball quotients; see \S\ref{sec:Hirz} for more on connections with his work. In Appendix A, we also describe why $\Gam_1$ is a normal subgroup of one of the lattices in $\PU(2,1)$ constructed by Deligne and Mostow \cite{DeligneMostow, MostowINT}. Specifically, the purpose of the appendix is to make this connection explicit enough that one can reproduce our Magma calculations using well-known presentations for Deligne--Mostow lattices.

\medskip

\medskip

We close with a question relating our work to the surfaces studied in Hirzebruch's paper \cite{Hirzebruch}.

\begin{qtn}\label{qtn:HirzebruchRigid}
The surfaces $X_n$ studied by Hirzebruch \cite{Hirzebruch} with the property that $c_1^2(X_n) / c_2(X_n) \to 3$ are $n^2$-fold \'etale covers of our surfaces $V_n$.
\begin{enum}

\item Are the surfaces $X_n$ also rigid?
\item Are the surfaces $V_n$ \'etale rigid?

\end{enum}
\end{qtn}

Briefly, a variety is \emph{\'etale rigid} if every \'etale cover is rigid. See \cite[Def.\ 2.1(9)]{BauerCatanese}. We note that the only known surfaces that are rigid but not \'etale rigid are Beauville surfaces. Also, notice that a negative answer to Question \ref{qtn:HirzebruchRigid}(1) immediately implies a negative answer for Question \ref{qtn:HirzebruchRigid}(2). Moreover, rigidity of $X_n$ implies rigidity of $V_n$ by an observation of Bauer--Catanese \cite[Prop.\ 2.5]{BauerCatanese}.

\subsubsection*{Acknowledgments}
The authors thank Fabrizio Catanese for conversations related to this paper.

\section{Hirzebruch's ball quotients}\label{sec:Hirz}

In this section we consider noncompact ball quotient manifolds first appearing, albeit not in this language, in work of Hirzebruch \cite{Hirzebruch}. We describe their explicit geometric construction, but we will also need presentations for the associated lattices in $\PU(2,1)$.

Let $\zeta = e^{\pi i / 3}$, $T$ denote the elliptic curve $\bbC / \bbZ[\zeta]$, and $A$ be the abelian surface $T \times T$ with coordinates $(z, w)$. Then $T_\al$ will denote the curve $w = \al z$ on $A$ and $T_\infty$ will denote the curve $z = 0$. Fix $n \ge 1$. If $U_n$ is the group of $n$-division points on $A$, then $D_\al^{(n)} = U_n(T_\al)$ is a divisor with support $n^2$ disjoint smooth irreducible curves.

Let $Y_n$ denote the blowup of $A$ at the $n^4$ points in $U_n$ and $E_j$ be the exceptional divisor above $j \in U_n$. If $\wt{D}_\al^{(n)}$ is the proper transform of $D_\al^{(n)}$ to $Y_n$, we define:
\[
Z_n = Y_n \ssm \bigcup_{\al \in \{0, 1, \zeta, \infty\}} \wt{D}_\al^{(n)}
\]
As noted in the acknowledgments to \cite{Hirzebruch}, $Z_n$ is a smooth, finite volume, noncompact quotient of the ball $\bbB^2$ by a torsion-free lattice $\Gam_n < \PU(2,1)$. Notice that $Z_n$ is an \'etale cover of $Z_1$ with covering group $(\bbZ / n)^4$, hence $\Gam_n$ is a normal subgroup of $\Gam_1$ with index $n^4$ and quotient $(\bbZ/n)^4$.

\medskip

In Appendix A, we describe how one can use Magma \cite{Magma} to show:

\begin{prop}\label{prop:HirzebruchLattice}
The lattice $\Gam_1 < \PU(2,1)$ has presentation with generators $h_1, \dots, h_4$ and relations:
\begin{align*}
h_3 h_2^{-1} h_1 h_4 h_2 h_1^{-1} h_3^{-1} h_4^{-1} &= \\
h_2 h_1^{-1} h_4^{-1} h_2^{-1} h_1 h_3 h_4 h_3^{-1} &= \\
h_2^{-1} h_1 h_3 h_1^{-1} h_3^{-1} h_4 h_2 h_4^{-1} &= \\
h_2 h_4 h_3^{-1} h_1^{-1} h_3 h_1 h_2^{-1} h_4^{-1} &= \\
h_1^{-1} h_3^{-1} h_2^{-1} h_4^{-1} h_3 h_4 h_1 h_2 &= \\ 
h_3^{-1} h_2^{-1} h_1 h_3 h_4^{-1} h_1^{-1} h_4 h_2 &= \\
h_1^{-1} h_3 h_1 h_3^{-1} h_4 h_2^{-1} h_4^{-1} h_2 &= \\
h_2^{-1} h_3 h_1 h_4^{-1} h_1^{-1} h_2 h_4 h_3^{-1} &= \\
h_3^{-1} h_4 h_2 h_4^{-1} h_3 h_4 h_3^{-1} h_2^{-1} h_3 h_4^{-1} &= \mathrm{Id}
\end{align*}
\end{prop}

Note that $\bbB^2 / \Gam_1$ has four cusps, each of which is smoothly compactified by an elliptic curve on $Y_n$ of self-intersection $-1$ that is isomorphic to $T$. See \cite{AMRT} for the basic theory of smooth toroidal compactifications. Consider the elements $k_1 := h_4 h_3^{-1}$, $k_2 := h_3 h_4^{-1}$, $z_1 := [k_1, k_2]$,
\begin{align*}
g_1 &:= k_1 k_2^{-1} & g_5 &:= x g_1 x^{-1} \\
g_2 &:= g_1 k_1^{-1} z_1^2 & g_6 &:= x g_2 x^{-1} \\
g_3 &:= x^{-1} g_1 x & g_7 &:= g_1^{-1} y g_1 y^{-1} g_1 \\
g_4 &:= x^{-1} g_2 x & g_8 &:= g_1^{-1} y g_2 y^{-1} g_1 \\
\\
w_1 &:= [g_1, g_2] & w_3 &:= [g_5, g_6] \\
w_2 &:= [g_3, g_4] & w_4 &:= [g_7, g_8]
\end{align*}
of $\Gam_1$, where $x,y$ are generators for the Deligne--Mostow lattice containing $\Gam_1$ considered in Appendix A. (Some of the above elements are very complicated as words in the generators $h_j$, so we leave it to the interested reader to calculate this in Magma.) A more detailed Magma analysis shows that the four conjugacy classes of parabolic subgroups of $\Gam_1$ associated with the cusps of $\bbB^2 / \Gam_1$ are generated by $g_{2j-1}$, $g_{2j}$, and $w_j$ subject to the relations
\[
[g_{2j-1}, g_{2j}] w_j^{-1} = [g_{2j-1}, w_j] = [g_{2j}, w_j] = \Id,
\]
$1 \le j \le 4$. Additionally, $w_1 w_2 w_3 w_4 = \Id$.

One can also see using Magma that any pair of distinct cusp subgroups generates $\Gam_1$. In particular, we can take $\{g_1, g_2, g_3, g_4\}$ as generators for $\Gam_1$. We also have:

\begin{lem}\label{lem:Commutator}
The commutator subgroup $\Gam_1^\prime$ of $\Gam_1$ is the normal subgroup $\langle\langle w_1, \dots, w_4 \rangle \rangle$ generated by the centers of its cusp subgroups.
\end{lem}

\begin{pf}
Killing the centers of the cusp subgroups of $\Gam_1$ is precisely the kernel of the induced map on fundamental groups from $\bbB^2 / \Gam_1$ to its smooth toroidal compactification $Y_1$. Since $\pi_1(Y_1) \cong \bbZ^4 \cong \Gam_1^{ab}$, the lemma follows.
\end{pf}

It is no accident that the abelianization of $\Gam_1$ is closely related to that of the fundamental group of its smooth toroidal compactification. In general, if $\bbB^2 / \Gam$ is a smooth ball quotient admitting a smooth toroidal compactification by the projective surface $V$, then $\bbB^2 / \Gam \hookrightarrow V$ induces a surjection on fundamental groups and an isomorphism $H_1(\Gam, \bbQ) \cong H_1(V, \bbQ)$; see \cite{DiCerboStoverCommutator} for an elementary proof. In other words, $\Gam^{ab}$ and $\pi_1(V)^{ab}$ always have the same free rank, though $\Gam^{ab}$ typically has larger torsion part.

\begin{lem}\label{lem:GenerateGam_n}
For all $n \ge 1$, the lattice $\Gam_n \trianglelefteq \Gam_1$ is generated by $\Gam_1^\prime$ and $\{g_1^n, \dots, g_4^n\}$.
\end{lem}

\begin{pf}
We will see that this follows from the fact that $\Gam_n$ is the kernel of the unique (up to choice of basis) homomorphism
\[
\rho_n : \Gam_1 \to (\bbZ / n)^4.
\]
Let $\rho_\infty : \Gam_1 \to \Gam_1^{ab} \cong \bbZ^4$ be the abelianization homomorphism. Recall from above that $g_1, \dots, g_4$ generate $\Gam_1$, hence $\rho_\infty(g_1), \dots, \rho_\infty(g_4)$ generate $\bbZ^4$. Then $\rho_\infty(g_1^n), \dots, \rho_\infty(g_4^n)$ are generators for the kernel of the homomorphism $\bbZ^4 \to (\bbZ / n)^4$ that sends each of $\rho_\infty(g_1), \dots, \rho_\infty(g_4)$ to one of the four standard generators for $(\bbZ / n)^4$.

Since $\Gam_n = \ker(\rho_n)$, for each $\gam \in \Gam_n$ we see that there are integers $j_1, \dots, j_4$ so that
\[
\rho_\infty(\gam) = \rho_\infty(g_1^n)^{j_1} \rho_\infty(g_2^n)^{j_2} \rho_\infty(g_3^n)^{j_3} \rho_\infty(g_4^n)^{j_4} \in \ker\left(\bbZ^4 \to (\bbZ / n)^4\right).
\]
Therefore, Lemma \ref{lem:Commutator} implies that
\[
\gam = (g_1^n)^{j_1} (g_2^n)^{j_2} (g_3^n)^{j_3} (g_4^n)^{j_4} w
\]
for some $w \in \Gam_1^\prime$. Thus $\Gam_n$ is contained in the subgroup of $\Gam_1$ generated by $g_1^n, \dots, g_4^n$ along with $\Gam_1^\prime$, but $g_1^n, \dots, g_4^n$ and $\Gam_1^\prime$ are all clearly contained in $\Gam_n$, so the proof of the lemma is complete.
\end{pf}

We now consider the structure of the cusp subgroups of $\bbB^2 / \Gam_n$ in relationship with the divisor $\calD_n := \sum_\al D_\al^{(n)}$. We start with the case $n = 1$. For each slope $\al \in \{0, 1, \zeta, \infty\}$, there is a unique $1 \le j(\al) \le 4$ so that the cusp subgroup of $\Gam_1$ associated with $T_\al$ is the group generated by $g_{2 j(\al) - 1}$ and $g_{2 j(\al)}$ with center generated by $w_{j(\al)}$. The correspondence between $\{0, 1, \zeta, \infty\}$ and $\{1, \dots, 4\}$ will not be of consequence for the arguments in this paper (in fact, the results in Appendix A imply that $\Aut(\bbB^2 / \Gam_1)$ acts transitively on the cusps, so the assignment is basically arbitrary). For $\bbB^2 / \Gam_n$, we have:

\begin{lem}\label{lem:CuspSubgroups}
The $4 n^2$ distinct conjugacy classes of cusp subgroups of $\Gam_n$ can be described as follows. For each $\al \in \{0, 1, \zeta, \infty\}$, there is a unique $j = j(\al) \in \{1, \dots, 4\}$ and there exist elements $\de_1, \dots, \de_{n^2}$ in $\Gam_1$ (depending on $j$) so that
\[
\left\{\de_i \langle g_{2j-1}^n, g_{2j}^n, w_j \rangle \de_i^{-1}\right\}_{i = 1}^{n^2}
\]
are representatives for the $\Gam_n$-conjugacy classes of cusp subgroups of $\Gam_n$ associated with the $n^2$ elliptic curves in the support of $D_\al^{(n)}$.
\end{lem}

\begin{pf}
This is an immediate consequence of elementary covering space theory applied to the regular cover $\bbB^2 / \Gam_n \to \bbB^2 / \Gam_1$ with group $(\bbZ / n)^4$. That the intersection of $\Gam_n$ with $\langle g_{2j-1}, g_{2j} \rangle$ is $\langle g_{2j-1}^n, g_{2j}^n, w_j \rangle$ follows directly from Lemma \ref{lem:GenerateGam_n}. To be precise, every element $\gam$ of the two-step nilpotent group $\langle g_{2j-1}, g_{2j} \rangle$, which is isomorphic to the usual $3$-dimensional integral Heisenberg group, has a unique normal form
\[
\gam = g_{2 j - 1}^{\ell_1} g_{2 j}^{\ell_2} w_j^{\ell_3}
\]
for $\ell_1, \ell_2, \ell_3 \in \bbZ$ (cf.\ the proof of Lemma \ref{lem:GenerateGam_n}). One can see this from either very general results on presentations of polycyclic groups (e.g., as in \cite[\S 9.4]{Sims}), or explicitly in this case from realizing $\langle g_{2j-1}, g_{2j} \rangle$ as the subgroup of upper-triangular matrices in $\SL_3(\bbZ)$. From this, we see that the kernel of restriction of the homomorphism $\rho_n$ in the proof of Lemma \ref{lem:GenerateGam_n} is precisely $\langle g_{2j-1}^n, g_{2j}^n, w_j \rangle$.
\end{pf}

Finally, we observe the following.

\begin{lem}\label{lem:Gam_n^ab}
One has $\Gam_n^{ab} \cong \bbZ^4 \oplus \calT$, where $\calT$ is a torsion abelian group of exponent dividing $n^2$.
\end{lem}

\begin{pf}
As mentioned above, the fact that $Z_n$ has smooth toroidal compactification birational to an abelian surface implies that $\Gam_n^{ab} / \calT \cong \bbZ^4$. By Lemma \ref{lem:GenerateGam_n} and Lemma \ref{lem:CuspSubgroups}, $\Gam_n$ is generated by $g_1^n, \dots, g_4^n$ along with elements of the form $\de_i w_j \de_i^{-1}$ for $\de_i \in \Gam_1$. Here we possibly introduce unnecessary generators and assume that we have a representative $\de_i w_j \de_i^{-1}$ associated with every $\Gam_n$-conjugacy class of cusp subgroups of $\Gam_n$.

The quotient of $\Gam_n$ by the elements $\{\de_i w_j \de_i^{-1}\}$ is then the fundamental group of the smooth toroidal compactification of $\bbB^2 / \Gam_n$, which is $\bbZ^4$. To prove the lemma it therefore suffices to show that each $\de_i w_j \de_i^{-1}$ has finite order dividing $n^2$ in $\Gam_n^{ab}$. This is clear from the easy calculation
\begin{equation}\label{eq:CuspComm}
(\de_i w_j \de_i^{-1})^{n^2} = \left[\de_i g_{2j-1}^n \de_i^{-1}, \de_i g_{2j}^n \de_i^{-1}\right]
\end{equation}
inside the two-step nilpotent group $\de_i \langle g_{2j-1}^n, g_{2j}^n, w_j \rangle \de_i^{-1} < \Gam_n$.
\end{pf}

\begin{rem}
The proof that $\de_i w_j \de_i^{-1}$ has order dividing $n^2$ in $\Gam_n^{ab}$ is closely related via Equation \eqref{eq:CuspComm} to the fact that the associated cusp is compactified by an elliptic curve of self-intersection $-n^2$.
\end{rem}

\section{Cyclic branched covers}\label{sec:Cyclic}

We freely use notation from \S\ref{sec:Hirz}. Consider an $n$-fold cyclic cover $V_n \to Y_n$ branched over the divisor $\calD_n := \sum_\al \wt{D}_\al^{(n)}$. If $W_n \subset V_n$ is the complement of the branch locus in $V_n$, then $W_n \to Z_n = \bbB^2 / \Gam_n$ is an \'etale cover. In other words, there is a normal subgroup $\Del_n \triangleleft \Gam_n$ with $\Gam_n / \Del_n \cong \bbZ / n$ so that $W_n = \bbB^2 / \Del_n$. An $n$-fold cyclic branched covering of $Y_n$ with branch set a (possibly empty) subset of $\calD_n$ is then equivalent to a certain epimorphism $\chi : \Gam_n \to \bbZ / n$. We call this the \emph{character} of the branched cover.

\medskip

The branched cover $V_n$ will be defined using a homomorphism $\tau$ from $\Gam_1$ to $\SL_5(\bbZ)$. Reduction modulo $n$ then defines a homomorphism $\tau_n$ from $\Gam_1$ to $\SL_5(\bbZ / n)$. When $n$ is odd, the kernel of $\tau_n$ will be an index $n$ subgroup of $\Gam_n$ that therefore defines an $n$-fold cyclic cover $V_n \to Y_n$.

Consider the representation $\tau : \Gam_1 \to \SL_5(\bbZ)$ defined by:
\begin{align*}
\tau(h_1) &= \begin{pmatrix} 1 & 0 & 1 & 1 & 1 \\ 0 & 1 & 0 & 0 & -1 \\ 0 & 0 & 1 & 0 & 1 \\ 0 & 0 & 0 & 1 & 1 \\ 0 & 0 & 0 & 0 & 1\end{pmatrix} & \tau(h_2) &= \begin{pmatrix} 1 & 1 & 0 & 0 & 0 \\ 0 & 1 & 0 & 0 & 0 \\ 0 & 0 & 1 & 0 & 0 \\ 0 & 0 & 0 & 1 & 0 \\ 0 & 0 & 0 & 0 & 1 \end{pmatrix} \\
\tau(h_3) &= \begin{pmatrix} 1 & 0 & 0 & 0 & 0 \\ 0 & 1 & 0 & 0 & 0 \\ 0 & 0 & 1 & 0 & 0 \\ 0 & 0 & 0 & 1 & 1 \\ 0 & 0 & 0 & 0 & 1 \end{pmatrix} & \tau(h_4) &= \begin{pmatrix} 1 & 1 & 0 & 1 & 1 \\ 0 & 1 & 0 & 0 & 1 \\ 0 & 0 & 1 & 0 & 0 \\ 0 & 0 & 0 & 1 & 1 \\ 0 & 0 & 0 & 0 & 1 \end{pmatrix}
\end{align*}
One checks directly using Proposition \ref{prop:HirzebruchLattice} that this is a representation. We similarly compute that:
\begin{align*}
\tau(g_1) &= \begin{pmatrix} 1 & 2 & -1 & 0 & 2 \\ 0 & 1 & 0 & 0 & 2 \\ 0 & 0 & 1 & 0 & -1 \\ 0 & 0 & 0 & 1 & -1 \\ 0 & 0 & 0 & 0 & 1 \end{pmatrix} & \tau(g_2) &= \begin{pmatrix} 1 & 1 & 0 & 1 & -1 \\ 0 & 1 & 0 & 0 & 1 \\ 0 & 0 & 1 & 0 & 0 \\ 0 & 0 & 0 & 1 & 0 \\ 0 & 0 & 0 & 0 & 1 \end{pmatrix} \\
\tau(g_3) &= \begin{pmatrix} 1 & 1 & -1 & -1 & 3 \\ 0 & 1 & 0 & 0 & 1 \\ 0 & 0 & 1 & 0 & -1 \\ 0 & 0 & 0 & 1 & -2 \\ 0 & 0 & 0 & 0 & 1 \end{pmatrix} & \tau(g_4) &= \begin{pmatrix} 1 & 0 & -1 & -1 & 1 \\ 0 & 1 & 0 & 0 & 1 \\ 0 & 0 & 1 & 0 & -1 \\ 0 & 0 & 0 & 1 & -2 \\ 0 & 0 & 0 & 0 & 1 \end{pmatrix}
\end{align*}
\[
\tau(w_1) = \begin{pmatrix} 1 & 0 & 0 & 0 & 1 \\ 0 & 1 & 0 & 0 & 0 \\ 0 & 0 & 1 & 0 & 0 \\ 0 & 0 & 0 & 1 & 0 \\ 0 & 0 & 0 & 0 & 1 \end{pmatrix}
\]
Reduction modulo $n$ then defines a homomorphism $\tau_n : \Gam_1 \to \SL_5(\bbZ / n)$. We now study some basic properties of $\tau$ and $\tau_n$. The first follows from basic matrix calculations:

\begin{lem}\label{lem:Ker1}
The kernel of $\tau$ contains the normal subgroup of $\Gam_1$ generated by $\{w_2 w_1^{-1}, w_4 w_3^{-1}, w_3 w_1\}$. Thus the same holds for $\tau_n$.
\end{lem}

The next fact we need is also clear from Lemma \ref{lem:Ker1} and the above expression for $\tau(w_1)$.

\begin{lem}\label{lem:Ker2}
For each $1 \le j \le 4$, $\tau_n(w_j)$ has order exactly $n$.
\end{lem}

We now see the first indication as to why we assume that $n$ is odd. We will comment more precisely later on why this assumption is natural for the methods of our proof.

\begin{lem}\label{lem:Ker3}
If $n$ is odd, then $\tau_n(g_j)$ has order exactly $n$ for every $j$ in $\{1,\dots, 4\}$.
\end{lem}

\begin{pf}
A simple induction shows that we have:
\begin{align*}
\tau(g_1)^m &= \begin{pmatrix} 1 & 2m & -m & 0 & \frac{m(5m-1)}{2} \\ 0 & 1 & 0 & 0 & 2m \\ 0 & 0 & 1 & 0 & -m \\ 0 & 0 & 0 & 1 & -m \\ 0 & 0 & 0 & 0 & 1 \end{pmatrix} \\
\tau(g_2)^m &= \begin{pmatrix} 1 & m & 0 & m & \frac{m(m-3)}{2} \\ 0 & 1 & 0 & 0 & m \\ 0 & 0 & 1 & 0 & 0 \\ 0 & 0 & 0 & 1 & 0 \\ 0 & 0 & 0 & 0 & 1 \end{pmatrix} \\
\tau(g_3)^m &= \begin{pmatrix} 1 & m & -m & -m & m(2m+1) \\ 0 & 1 & 0 & 0 & m \\ 0 & 0 & 1 & 0 & -m \\ 0 & 0 & 0 & 1 & -2m \\ 0 & 0 & 0 & 0 & 1 \end{pmatrix} \\
\tau(g_4)^m &= \begin{pmatrix} 1 & 0 & -m & -m & \frac{m(3m-1)}{2} \\ 0 & 1 & 0 & 0 & m \\ 0 & 0 & 1 & 0 & -m \\ 0 & 0 & 0 & 1 & -2m \\ 0 & 0 & 0 & 0 & 1 \end{pmatrix}
\end{align*}
for any integer $m$. Then, $\tau_n(g_j)^m$ is clearly nontrivial in $\SL_5(\bbZ / n)$ for every $1 \le m < n$. When $n$ is odd, we also have that $\tau_n(g_j)^n$ is the identity. Indeed, the $(1,5)$ coordinate is divisible by $n$ if and only if $n$ is odd, and the other coordinates are clearly the same as the identity matrix. This proves the lemma.
\end{pf}

Let $G_n := \tau_n(\Gam_1) < \SL_5(\bbZ / n)$ and $\Del_n := \ker(\tau_n)$. We have:

\begin{lem}\label{lem:rho_nImage}
If $n$ is odd, there is a central exact sequence:
\[
1 \to \bbZ / n = \Gam_n / \Del_n \to G_n = \Gam_1 / \Del_n \to (\bbZ / n)^4 = \Gam_1 / \Gam_n \to 1
\]
In particular, $G_n$ is a two-step nilpotent group of order $n^5$ with abelianization $(\bbZ / n)^4$ and center $\bbZ / n$.
\end{lem}

\begin{pf}
First, notice that $\tau_n(w_j)$ is equal to either $\tau_n(w_1)$ or its inverse for each $1 \le j \le 4$ by Lemma \ref{lem:Ker1}. Since $w_j$ commutes with $g_{2 j - 1}$ and $g_{2 j}$, we see that $\tau_n(w_1)$ is central in $G_n$. It has order $n$ by Lemma \ref{lem:Ker2}. This, along with Lemma \ref{lem:Commutator}, implies that the image of $\Gam_1^\prime$ in $G_n$ is the cyclic subgroup generated by $\tau_n(w_1)$. The quotient group is $(\bbZ / n)^4$ by Lemma \ref{lem:Ker3}. That $G_n$ is a two-step nilpotent group is immediate.
\end{pf}

\begin{rem}\label{rem:NormalForm}
Since $G_n$ is a two-step nilpotent group with commutator subgroup generated by $\tau_n(w_1)$, every element can be put in the normal form
\[
\tau_n(g_1)^{m_1} \tau_n(g_2)^{m_2} \tau_n(g_3)^{m_3} \tau_n(g_4)^{m_4} \tau_n(w_1)^{m_5}
\]
for $0 \le m_j < n$. As a matrix in $\SL_5(\bbZ / n)$, this element equals
\[
\left(\begin{smallmatrix} 1 & 2 m_1 + m_2 + m_3 & -m_1 - m_3 - m_4 & m_2 - m_3 - m_4 & t_n \\
0 & 1 & 0 & 0 & 2 m_1 + m_2 + m_3 + m_4 \\
0 & 0 & 1 & 0 & -m_1 - m_3 - m_4 \\
0 & 0 & 0 & 1 & -m_1 - 2 m_3 - 2 m_4 \\
0 & 0 & 0 & 0 & 1 \end{smallmatrix}\right)
\]
where
\begin{align*}
t_n := \frac{1}{2} (-m_1 &+ 5 m_1^2 - 3 m_2 + 4 m_1 m_2 + m_2^2 \\
&+ 2 m_3 + 6 m_1 m_3 - 2 m_2 m_3 + 4 m_3^2 - m_4 \\
&+ 6 m_1 m_4 - 2 m_2 m_4 + 8 m_3 m_4 + 3 m_4^2 + 2 m_5).
\end{align*}
Note that $n$ is odd, so $2$ is invertible. It is not hard to deduce from this matrix representation that each choice of $m_1, \dots, m_5$ gives a unique element of $G_n$, hence $G_n$ has order $n^5$.
\end{rem}

We are now prepared to connect $\tau_n$ to an $n$-fold cyclic cover $V_n$ of $Y_n$.

\begin{prop}\label{prop:OurCover}
For $n \ge 2$ odd, let $\chi : \Gam_n \to \bbZ / n$ be the restriction of $\tau_n$ to $\Gam_n$. Then $\chi$ is the character of an $n$-fold cyclic branched cover $V_n \to Y_n$ branched over $\calD_n$ that can be described as follows:
\begin{enum}

\item There is a generator $\sig$ for $\bbZ / n$ so that for each $\al \in \{0, 1, \zeta, \infty\}$, we have $\chi(w_{j(\al)}) = \sig$ for $j(\al) = 1,2$ and $\chi(w_{j(\al)}) = \sig^{-1}$ for $j(\al)=3,4$.

\item For each $\al \in \{0, 1, \zeta, \infty\}$, the cover has ``the same character'' over each elliptic curve in the support of $D_\al^{(n)}$. That is, for each cusp subgroup $\{\de \langle g_{2j(\al)-1}^n, g_{2j(\al)}^n, w_{j(\al)}\rangle \de^{-1}\}$ of $\Gam_n$ associated with one of the elliptic curves in the support of $D_\al^{(n)}$ as in Lemma \ref{lem:CuspSubgroups}, we have that $\chi(\de w_{j(\al)} \de^{-1}) = \chi(w_{j(\al)})$.

\end{enum}
Consequently, $V_n \to Y_n$ is branched to order precisely $n$ over each of the elliptic curves in the support of $\calD_n$. Moreover, $V_n$ is the smooth toroidal compactification of the ball quotient $\bbB^2 / \Del_n$.
\end{prop}

\begin{pf}
From Lemma \ref{lem:rho_nImage} we have that $\chi$ defines a homomorphism from $\Gam_n$ to $\bbZ / n$, hence it describes an $n$-fold cyclic cover $V_n$ of $Y_n$. Set $\sig := \tau_n(w_1)$. By Lemma \ref{lem:Ker1} we have that
\[
\chi(w_{j(\al)}) = \begin{cases} \sig & j(\al) = 1,2 \\ \sig^{-1} & j(\al) = 3,4 \end{cases}
\]
which gives the first part of the description of $V_n$. The second part follows immediately from the fact that $\de w_{j(\al)} \de^{-1}$ describes a small loop around the appropriate elliptic curve in the support of $D_\al^{(n)}$ and that
\[
\tau_n(\de w_{j(\al)} \de^{-1}) = \tau_n(w_{j(\al)})
\]
for any $\de \in \Gam_1$, since $\tau_n(w_{j(\al)})$ is central in $G_n$.

The claim about the branching of $V_n \to Y_n$ is immediate from this description. That $V_n$ is the smooth toroidal compactification of $\bbB^2 / \Del_n$ is an immediate consequence of our construction, since $\bbB^2 / \Gam_n$ is precisely the complement in $Y_n$ of the branch set. This proves the proposition.
\end{pf}

Our goal now is to compute the irregularity of $V_n$. The key technical result is the following.

\begin{prop}\label{prop:GenerateDeln}
Suppose $n$ is odd. The kernel $\Del_n$ of $\tau_n$ is generated by $g_1^n, \dots, g_4^n$ along with the normal closure in $\Gam_1$ of
\[
\calW := \{w_1^n, \dots, w_4^n, w_2 w_1^{-1}, w_4 w_3^{-1}, w_3 w_1\}.
\]
\end{prop}

\begin{pf}
Let $\Lam_n$ be the subgroup of $\Gam_1$ generated by $g_1^n, \dots, g_4^n$ along with the normal closure in $\Gam_1$ of $\calW$. Lemmas \ref{lem:Ker1} - \ref{lem:Ker3} imply that $\Lam_n < \Del_n$. In fact, they imply that $\Del_n$ contains the normal closure of $\Lam_n$ in $\Gam_1$.

Consider the group
\[
H := \Gam_1 / \langle \langle w_2 w_1^{-1}, w_4 w_3^{-1}, w_3 w_1 \rangle \rangle.
\]
Note that $\langle \langle w_2 w_1^{-1}, w_4 w_3^{-1}, w_3 w_1 \rangle \rangle \le \Lam_n$ by definition. Moreover, note that the quotient homomorphism $\Gam_1 \to G_n := \Gam_1 / \Del_n$ factors through $\Gam_1 \to H$ by construction.

First, we note that $H$ is a two-step nilpotent group. Let $\wh{g}_j$ be the image of $g_j$ in $H$ and $\wh{w}$ the image of $w_1$. The relations defining $H$ imply that each $\wh{g}_j$ commutes with $\wh{w}$, since each $\wh{w}_\ell$ maps to either $\wh{w}$ or its inverse, and $\wh{w}_\ell$ commutes with the generators $\wh{g}_{2 \ell - 1}$ and $\wh{g}_{2 \ell}$. Then $\wh{w}$ is central in $H$ with $H / \langle \langle \wh{w} \rangle \rangle \cong \bbZ^4$ by Lemma \ref{lem:Commutator}, hence $H$ is two-step nilpotent.

Since $H$ is a two-step nilpotent group generated by $\wh{g}_1, \dots, \wh{g}_4$ with center generated by $\wh{w}$, every $h \in H$ has a unique expression in normal form as
\[
h = \wh{g}_1^{j_1} \wh{g}_2^{j_2} \wh{g}_3^{j_3} \wh{g}_4^{j_4} \wh{w}^k
\]
for $j_1, \dots, j_4$ and $k$ all integers. This follows from very general results on normal form for torsion-free polycyclic groups; e.g., see \cite[\S 9.4]{Sims}. One can also see this by mapping $H$ into $\SL_5(\bbZ)$ under $\tau$, then using the matrix formula in Remark \ref{rem:NormalForm} to show that $H \cong \tau(\Gam_1)$.

The kernel of the homomorphism from $H$ onto $G_n$ is clearly the image of $\Del_n$ in $H$. This kernel is also given by setting $\wh{w}$ and each $\wh{g}_j$ to have order $n$, that is, it is generated by $\wh{g}_1^n, \dots, \wh{g}_4^n$ and $\wh{w}^n$ (i.e., it is generated by those five elements, not their normal closure). This is precisely the image of $\Lam_n$ in $H$. Since $\ker(\Gam_1 \to H) \le \Lam_n \le \Del_n$, the usual subgroup correspondence for quotient groups implies that $\Lam_n = \Del_n$, which proves the proposition.
\end{pf}

We now arrive at one of the key results about $\Del_n$.

\begin{prop}\label{prop:ComputeAb}
One has that $\Del_n^{ab} / \mathrm{Torsion} \cong \bbZ^4$.
\end{prop}

\begin{pf}
By Proposition \ref{prop:GenerateDeln}, we have that $\Del_n^{ab}$ is generated by the images of the four elements $g_1^n, \dots, g_4^n$ along with elements of the form
\[
\de w_j^n \de^{-1},\, \de w_2 w_1^{-1} \de^{-1},\, \de w_4 w_3^{-1} \de^{-1},\, \textrm{or}\ \de w_3 w_1 \de^{-1}
\]
for some $\de \in \Gam_1$. Since $g_1, \dots, g_4$ generate $\Gam_1^{ab} \cong \bbZ^4$, we see that $g_1^n, \dots, g_4^n$ generate a $\bbZ^4$ in $\Del_n^{ab}$. To prove the proposition, it then suffices to show that all the other generators have finite order in $\Del_n^{ab}$.

As in the proof of Lemma \ref{lem:Gam_n^ab}, the appropriate power of any such element is a product of commutators of elements of $\Del_n^{ab}$. Indeed, we have
\begin{align*}
\left[\de g_{2 j - 1}^n \de^{-1}, \de g_{2 j}^n \de^{-1} \right] &= \left(\de w_j \de^{-1}\right)^{n^2} \\
&= \left(\de w_j^n \de^{-1}\right)^n,
\end{align*}
so $\de w_j^n \de^{-1}$ has order divisible by $n$ in $\Del_n^{ab}$. Then $(\de w_3 w_1^{-1} \de^{-1})^n$ has the same image in $\Del_n^{ab}$ as $(\de w_3^n \de^{-1})(\de w_1^{-n} \de^{-1})$, so $\de w_3 w_1^{-1} \de^{-1}$ has order divisible by $n^2$ in $\Del_n^{ab}$, and the same argument works for $\de w_4 w_2^{-1} \de^{-1}$ and $\de w_2 w_1 \de^{-1}$. This proves the corollary.
\end{pf}

This immediately implies the main result of this section:

\begin{cor}\label{cor:Irregularity}
The algebraic surface $V_n$ has irregularity $2$.
\end{cor}

\begin{pf}
This follows from Proposition \ref{prop:ComputeAb} and the fact, referenced above, that $H_1(\Del_n, \bbQ)$ is isomorphic to $H_1(V_n, \bbQ)$.
\end{pf}

\begin{rem}
It is absolutely critical to our argument that $g_1^n, \dots g_4^n$ be generators, not just normal generators, for $\Del_n$ (along with the normal closure of the set $\calW$). Otherwise, the proof of Proposition \ref{prop:ComputeAb} would completely fall apart. If we only knew that $\Del_n$ were normally generated by the $g_j$s, then distinct generators $g_j^n$ and $\de g_j^n \de^{-1}$ could determine linearly independent elements of $\Del_n^{ab}$ of infinite order, hence our conclusion that $\Del_n^{ab}$ modulo the images of the normal closure of $\calW$ is isomorphic to $\bbZ^4$ would be incorrect.
\end{rem}

We also must prove:

\begin{prop}\label{prop:ChernCalc}
The surfaces $V_n$ defined above are distinct minimal surfaces of general type for which $\frac{c_1^2(V_n)}{c_2(V_n)} = 3-\frac{4}{n^2}$.
\end{prop}

\begin{pf}
This follows analogously to calculations done by Hirzebruch in \cite[\S 1]{Hirzebruch}. We first argue $V_n$ is minimal by showing that it contains no curves of genus zero. The only curves of genus zero on $Y_n$ are the exceptional curves for the blowup $Y_n \to A$, where $A$ is the abelian surface from \S \ref{sec:Hirz}. The Hurwitz formula gives that each exceptional curve on $V_n$ lifts to a curve of genus $n-1$ on $V_n$, since the branch locus meets each exceptional curve in exactly four points. Since $n \ge 3$, we see that $V_n$ is minimal.

Recalling notation from \S\ref{sec:Hirz}, we have:
\begin{align*}
c_1^2(V_n) &= n \left(\sum_{j \in U_n} E_j + \left(1 - \frac{1}{n}\right) \sum_{\al \in \{0, 1, \zeta, \infty\}} \wt{D}_\al^{(n)}\right)^2 \\
&= 3 n^5 - 4 n^3 \\
c_2(V_n) &= n^5
\end{align*}
In the calculation of $c_1^2$ we used that each $E_j$ is an exceptional curve, $\sum \wt{D}_\al^{(n)}$ consists of $4 n^2$ disjoint elliptic curves of self-intersection $-n^2$, and that each curve in the support of $\wt{D}_\al^{(n)}$ intersects exactly $n^2$ of the exceptional curves. The computation of $c_2(V_n)$ is easier, since $c_2(Y_n) = n^4$ and $V_n \to Y_n$ is branched over elliptic curves.

Since $V_n$ is minimal and contains no curves of genus zero, the fact that $c_1^2(V_n)$ and $c_2(V_n)$ are positive implies that $V_n$ is of general type. They are clearly distinct surfaces, as their Chern numbers are different for each $n$. This completes the proof.
\end{pf}




We now prove the last of the stated properties of $V_n$ other than rigidity.

\begin{lem}\label{lem:AmpleK}
The canonical divisor on $V_n$ is ample.
\end{lem}

\begin{pf}
Since $V_n$ is minimal of general type, we know that $K_{V_n}^2 > 0$. The Nakai--Moishezon criterion implies that we need only show that $K_{V_n}$ meets every irreducible curve on $V_n$ positively. Since $\sum D_\al^{(n)}$ is ample on the abelian variety $A$, it moreover suffices to notice that
\[
\left(\sum_{j \in U_n} E_j + \left(1 - \frac{1}{n}\right) \sum_{\al \in \{0, 1, \zeta, \infty\}} \wt{D}_\al^{(n)}\right)\cdot E_i = -1 + 4\left(1 - \frac{1}{n}\right) > 0
\]
for every exceptional curve $E_i$ on $Y_n$.
\end{pf}

\begin{rem}
It seems to be the case that $V_n$ has contractible universal cover for all odd $n$, but we have not checked this carefully. If true, then $V_n$ is a projective classifying space in the sense of \cite[Def.\ 3.2]{BauerCatanese}. Bauer and Catanese also note that all known examples of rigid surfaces are projective classifying spaces. For example, note that $V_n$ does not contain rational curves, so the obvious obstruction to a projective surface having nontrivial $\pi_2$ fails.

One can see that $V_n$ is aspherical for sufficiently large $n$ using work of Hummel and Schroeder \cite{HummelSchroeder}. In particular, the cusp subgroups of $\Del_n$ are generated by conjugates of $w_j^n$, which for $n$ sufficiently large will have $|w_j^n| > \rho$ in the notation of \cite[\S3]{HummelSchroeder}. Following \cite[Remark 1]{HummelSchroeder}, the smooth toroidal compactification of $\bbB^2 / \Del_n$, which is precisely $V_n$, will admit a metric of nonpositive curvature. Then the Cartan--Hadamard theorem implies that the universal covering of $V_n$ is contractible, as desired.
\end{rem}

\begin{rem}
We close with a remark on the case where $n$ is even. We saw in the proof of Lemma \ref{lem:Ker3} that $\tau_n(g_j)$ does not have order $n$ when $n$ is even. This does not seem terribly problematic: We could simply define $\Del_n$ to be the subgroup of $\Gam_1$ generated by $g_1^n, \dots, g_4^n$ along with the normal closure of the set $\calW$ from the statement of Proposition \ref{prop:GenerateDeln} and not worry about $\tau_n$ (this was our original approach to the problem). Come to find out, $\Del_n$ is a normal subgroup of $\Gam_1$ but only of index $n^5 / 2$, since the image of $w_1$ in $\Gam_1 / \Del_n$ only has order $n / 2$. Consequently, one could perhaps replace $V_n$ with the $(n/2)$-fold cyclic cover branched over $\calD_n$ and the results in this section would carry through, but we did not explore this carefully.
\end{rem}

\section{Rigidity of the surface $V_n$}\label{sec:Rigidity}

In \S\ref{sec:Cyclic}, for every odd $n \ge 3$ we constructed a sequence $\{V_n\}$ of distinct minimal smooth projective surfaces of general type with irregularity two so that their universal cover is contractible and $\frac{c_1^2(V_n)}{c_2(V_n)} = 3-\frac{4}{n^2}$. To prove Theorem \ref{thm:Main}, it remains to show that $V_n$ is rigid. In this section we will prove that $V_n$ having irregularity two implies $H^1(V_n,T_{V_n})=0$, where $T_{V_n}$ is the tangent bundle of $V_n$, hence the surfaces $V_n$ are rigid. For a line bundle $\mathcal F$, we will write $\mathcal F^i := \mathcal F^{\otimes i}$.

We first algebraically construct a surface $S_n$ that is isomorphic to $V_n$. Let $E$ be the sum of the $n^4$ exceptional curves $E_j$ of the blowup $\sigma_n \colon Y_n \to A$. We have the linear equivalence
\[
\wt{D}_\al^{(n)} + E \sim n^2 \sigma_n^*(T_\al)
\]
for each $\al \in \{0,1,\infty, \zeta \}$, and therefore if
\[
\L := \O_{Y_n}\big(n\sigma_n^*(T_0)+n\sigma_n^*(T_1)+n(n-1)\sigma_n^*(T_\infty)+ n(n-1)\sigma_n^*(T_\zeta) -2E \big)
\]
we then have the isomorphism of line bundles
\[
\O_{Y_n}\big(\wt{D}_0^{(n)} + \wt{D}_1^{(n)} + (n-1) \wt{D}_\infty^{(n)} + (n-1) \wt{D}_\zeta^{(n)} \big)\simeq \L^{n}.
\]
With this data, we construct as in \cite[\S 3]{EsnaultViehweg} (see also \cite[Ch.\ IV]{Urzua}) a nonnormal irreducible projective $n$-th root cover $S'_n$ of $Y_n$ branched along $\wt{D}_0^{(n)} + \wt{D}_1^{(n)} + (n-1) \wt{D}_\infty^{(n)} + (n-1) \wt{D}_\zeta^{(n)}$.

The variety $S'_n$ is defined as $\spec_{\O_{Y_n}} \big(\oplus_{i=0}^{n-1} {\L}^{-i} \big)$, and the corresponding finite morphism $f'_n \colon S'_n \to Y_n$ gives ${f'_n}_* \O_{S'_n} \simeq \oplus_{i=0}^{n-1} \L^{-i}$. As in \cite[\S 3]{EsnaultViehweg}, the normalization $S_n \to S'_n$ of $S'_n$ is constructed as $S_n:= \spec_{\O_{Y_n}} \big(\oplus_{i=0}^{n-1} {\L^{(i)}}^{-1} \big)$ where
\[
\L^{(i)} := \L^{i} \otimes \O_{Y_n}\big(-(i-1)\wt{D}_\infty^{(n)} - (i-1) \wt{D}_\zeta^{(n)} \big),
\]
and the composition $f_n \colon S_n \to Y_n$ gives ${f_n}_* \O_{S_n} \simeq \oplus_{i=0}^{n-1} {\L^{(i)}}^{-1}$.

Note that the construction of $S_n$ does not depend on the assumption that $n$ is odd. When this is the case, the choice of multiplicities $1,1,n-1,n-1$ for the branch divisor $\wt{D}_0^{(n)} + \wt{D}_1^{(n)} + (n-1) \wt{D}_\infty^{(n)} + (n-1) \wt{D}_\zeta^{(n)}$ coincides with the choice of the character $\Gamma_n \to \bbZ/n$ in Proposition \ref{prop:OurCover}, after choosing the right bijection $j$. Therefore $S_n$ and $V_n$ are isomorphic. From now on we will use  $V_n$ instead of $S_n$ to keep the notation from previous sections.

For a nonsingular projective surface $S$ and a simple normal crossing divisor $D$, let $\Omega_{S}^1(\log D)$ be the sheaf of log differentials as in \cite[\S 2]{EsnaultViehweg}. By \cite[Thm.\ VI.5]{Urzua} (see also \cite[Lem.\ 3.1.4]{Zuo}, \cite[Prop.\ 4.1]{Pardini}, \cite[Prop.\ 5.7]{BauerCatanese}), we have
\[
{f_n}_* \big(\Omega_{V_n}^1 \otimes \Omega_{V_n}^2 \big) \simeq \bigoplus_{i=0}^{n-1} \big( \Omega_{Y_n}^1(\log D^{(i)}) \otimes \Omega_{Y_n}^2 \otimes \L^{(i)} \big)
\]
where $D^{(i)}= \wt{D}_0^{(n)} + \wt{D}_1^{(n)} + \wt{D}_\infty^{(n)} + \wt{D}_\zeta^{(n)}$ for $i\neq 1,n-1$, $D^{(1)}=\wt{D}_0^{(n)} + \wt{D}_1^{(n)}$, and $D^{(n-1)}=\wt{D}_\infty^{(n)} + \wt{D}_\zeta^{(n)}$. Thus, by Serre duality and since $f_n$ is affine, we have
\[
H^1(V_n,T_{V_n}) \simeq \bigoplus_{i=0}^{n-1} H^1(Y_n,\Omega_{Y_n}^1(\log D^{(i)}) \otimes \Omega_{Y_n}^2 \otimes \L^{(i)} ).
\]
We now show:

\begin{lem}\label{lem:tech1}
If we have that $H^1(Y_n,\Omega_{Y_n}^1 \otimes \Omega_{Y_n}^2\otimes {\L^{(i)}})=0$ for all $i \neq 0$, then $H^1(V_n,T_{V_n})=0$.
\end{lem}

\begin{pf}
For $i \neq 0$, we tensor the residue short exact sequence (see \cite[\S 2]{EsnaultViehweg})
\[
0 \to \Omega_{Y_n}^1 \to \Omega_{Y_n}^1(\log D^{(i)}) \to \oplus_{C \in D^{(i)}} \O_{C} \to 0
\]
with $\Omega_{Y_n}^2 \otimes {\L^{(i)}}$, where $C$ runs over the irreducible curves in $D^{(i)}$. Since $C^2=-n^2$ and $\text{deg}_{C} {\L^{(i)}}|_{C}$ is either $-in$ or $-n(n-i)$, we have that $\text{deg}_{C} \Omega_{Y_n}^2 \otimes {\L^{(i)}}|_{C}$ is either $n(n-i)>0$ or $in>0$, and hence $H^1(C, \Omega_{Y_n}^2 \otimes {\L^{(i)}}|_{C})=0$.

Therefore, by considering the long exact sequence in cohomology, we have that $H^1(Y_n,\Omega_{Y_n}^1 \otimes \Omega_{Y_n}^2\otimes {\L^{(i)}})=0$ implies that
\[
H^1(Y_n,\Omega_{Y_n}^1(\log D^{(i)}) \otimes \Omega_{Y_n}^2 \otimes \L^{(i)})=0.
\]
Now by Fujiki's Theorem \cite[Thm.\ 4.1]{Fujiki}, we have that
\[
H^1(Y_n,\Omega_{Y_n}^1(\log D^{(0)}) \otimes \Omega_{Y_n}^2)=0,
\]
since $Y_n \ssm D^{(0)}$ is a ball quotient (see \S\ref{sec:Cyclic}). As
\[
H^1(V_n,T_{V_n}) \simeq \bigoplus_{i=0}^{n-1} H^1(Y_n,\Omega_{Y_n}^1(\log D^{(i)}) \otimes \Omega_{Y_n}^2 \otimes \L^{(i)} ),
\]
we finally obtain $H^1(V_n,T_{V_n})=0$.
\end{pf}

\begin{lem}\label{lem:tech2}
Given $i \neq 0$, we have that $H^1(Y_n, {\L^{(i)}}^{-1})=0$ if and only if $H^1(Y_n,\Omega_{Y_n}^1 \otimes \Omega_{Y_n}^2 \otimes {\L^{(i)}})=0$.
\end{lem}

\begin{pf}
We have that $\sigma_n^*(\Omega_{T\times T}^1)$ can be identified in $\Omega_{Y_n}^1$ with $\Omega_{Y_n}^1(\log E) \otimes \O_{Y_n}(-E)$ (e.g., see \cite[Lem.\ 3.1.3]{Zuo}), and so via the short exact sequence
\[
0 \to \Omega_{Y_n}^1(\log E) \otimes \O_{Y_n}(-E) \to \Omega_{Y_n}^1 \to \bigoplus_{j=1}^{n^4} \Omega_{E_j}^1 \to 0
\]
(see \cite[\S 2]{EsnaultViehweg}) we obtain the short exact sequence
\[
0 \to \sigma_n^*(\Omega_{T \times T}^1) \to \Omega_{Y_n}^1 \to \bigoplus_{j=1}^{n^4} \Omega_{E_j}^1 \to 0.
\]
However, $\Omega_{T \times T}^1 = \O_{T \times T} \oplus \O_{T \times T}$, and so it becomes
\[
0 \to \O_{Y_n} \oplus \O_{Y_n} \to \Omega_{Y_n}^1 \to \bigoplus_{j=1}^{n^4} \Omega_{E_j}^1 \to 0 .
\]
We now obtain
\begin{align*}
0 \to (\Omega_{Y_n}^2 \otimes {\L^{(i)}}) \oplus (\Omega_{Y_n}^2 \otimes {\L^{(i)}}) &\to \Omega_{Y_n}^1 \otimes \Omega_{Y_n}^2 \otimes {\L^{(i)}} \\
&\to \bigoplus_{j=1}^{n^4} (\Omega_{E_j}^1 \otimes \Omega_{Y_n}^2 \otimes {\L^{(i)}}) \to 0
\end{align*}
after tensoring with $\otimes (\Omega_{Y_n}^2 \otimes {\L^{(i)}})$. The associated long exact sequence then gives
\begin{align*}
\cdots \to H^1(Y_n,\Omega_{Y_n}^2 \otimes {\L^{(i)}})^2 &\to H^1(Y_n, \Omega_{Y_n}^1 \otimes \Omega_{Y_n}^2 \otimes {\L^{(i)}}) \to \\  &\to \bigoplus_{j=1}^{n^4} H^1(E_j, \Omega_{E_j}^1 \otimes \Omega_{Y_n}^2 \otimes {\L^{(i)}}) \to \cdots
\end{align*}
Since the multiplicities of the branch divisor are $1,1,n-1,n-1$, we have that the degree of $\Omega_{E_j}^1 \otimes \Omega_{Y_n}^2 \otimes {\L^{(i)}}$ in $E_j$ is equal to $-1$. Therefore, by Serre duality we obtain that $H^1(Y_n, {\L^{(i)}}^{-1})=0$ if and only if $H^1(Y_n,\Omega_{Y_n}^1 \otimes \Omega_{Y_n}^2 \otimes {\L^{(i)}})=0$, which completes the proof.
\end{pf}

\begin{thm}\label{lem:rigid}
The surface $V_n$ is rigid for every odd $n \ge 3$.
\end{thm}

\begin{pf}
Since ${f_n}_* \O_{V_n} \simeq \oplus_{i=0}^{n-1} {\L^{(i)}}^{-1}$ and $f_n$ is affine, we have
\[
H^1(V_n,\O_{V_n}) \cong \bigoplus_{i=1}^{n-1} H^1(Y_n, {\L^{(i)}}^{-1}).
\]
By Corollary \ref{cor:Irregularity}, we have that $\text{dim}_{\bbC} H^1(V_n,\O_{V_n})=2$. Since
\[
\text{dim}_{\bbC} H^1(Y_n,\O_{Y_n})=2,
\]
we must have $H^1(Y_n, {\L^{(i)}}^{-1})=0$ for all $i \neq 0$. Then, by Lemma \ref{lem:tech2} and Lemma \ref{lem:tech1}, this implies that $H^1(V_n,T_{V_n})=0$, and hence $V_n$ is a rigid surface.
\end{pf}

\begin{rem}
We point out that in the proof of rigidity we could have used Fujiki's theorem for all $i$ except $1$ and $n-1$, where it was essential to know that the irregularity of $V_n$ is equal to $2$.
\end{rem}

\section*{Appendix A: A Deligne--Mostow orbifold}

\renewcommand{\thethm}{\Alph{section}.\arabic{thm}}
\setcounter {section}{1}

In this appendix, we consider the Deligne--Mostow orbifold with associated weights $(5, 4, 1, 1, 1)$ \cite{DeligneMostow, MostowINT} and explicitly relate the ball quotient $\bbB^2 / \Gam_1$ from \S\ref{sec:Hirz} to finite coverings of this orbifold. One can use this to quickly deduce a presentation for $\Gam_1$ using computer algebra software like Magma \cite{Magma} and then verify the other group-theoretic claims made in this paper.

Nothing in this section is new, though we do not have a single reference for everything we use. We broadly follow the treatment of the structure of these orbifolds given by Kirwan, Lee, and Weintraub \cite{KirwanLeeWeintraub}. The weights $(5,4,1,1,1)$ satisfy Mostow's condition $\Sig$INT \cite{DeligneMostow, MostowINT}, and the underlying space for the orbifold is the quotient of $\bbP^2$ by the action of the symmetric group $S_3$ by coordinate permutations. This quotient is the weighted projective space $\bbP(1,2,3)$.

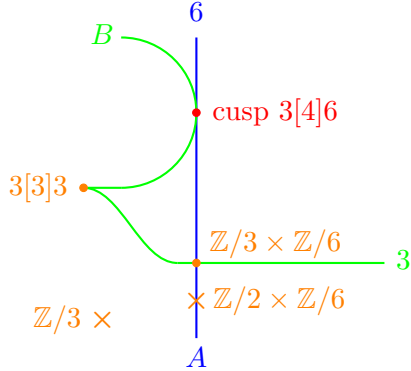
\begin{figure} [htbp]
\centering
\begin{tikzpicture}
\draw[thick, color=blue] (0,2) -- (0,-2);
\draw[thick, color=green] (-1,2) arc (90:-90:1cm);
\draw[thick, color=green] (-1, 0) -- (-1.5, 0);
\draw[thick, color=green] (-1.5,0) .. controls (-1, 0) and (-0.75, -1) .. (-0.25,-1);
\draw[thick, color=green] (-0.25, -1) -- (2.5,-1);
\node [color=blue] at (0,-2.25) {$A$};
\node [color=blue] at (0, 2.35) {$6$};
\node [color=green] at (-1.25, 2.05) {$B$};
\node [color=green] at (2.75, -0.95) {$3$};
\draw [color=red, fill=red] (0,1) circle (0.05cm);
\node [color=red] at (1.05, 1) {$\textrm{cusp}\ 3[4]6$};
\draw [color=orange, fill=orange] (-1.5,0) circle (0.05cm);
\node [color=orange] at (-2.1, 0) {$3[3]3$};
\draw [color=orange, fill=orange] (0,-1) circle (0.05cm);
\node [color=orange] at (1.05, -0.75) {$\bbZ / 3 \times \bbZ / 6$};
\node [color=orange] at (0, -1.5) {$\boldsymbol{\times}$};
\node [color=orange] at (1.1, -1.5) {$\bbZ / 2 \times \bbZ / 6$};
\node [color=orange] at (-1.25, -1.75) {$\boldsymbol{\times}$};
\node [color=orange] at (-1.85, -1.75) {$\bbZ / 3$};
\end{tikzpicture}
\caption{The orbifold structure on $\bbP(1,2,3)$}\label{fig:DM}
\end{figure}

We now describe the orbifold structure on $\bbP(1,2,3)$ associated with these weights. See Figure \ref{fig:DM}. Let $[x : y : z]$ denote homogeneous coordinates on $\bbP^2$. The curves $x=0$, $y=0$, and $z=0$ project to the curve $A$ on $\bbP(1,2,3)$ with an $\mathrm{A}_1$ singularity marked by a $\times$. Similarly, the curves $x=y$, $y=z$, and $z=x$ project to the curve $B$ with a cusp singularity that is a smooth point of the surface. Then $A$ has orbifold weight $6$ and $B$ has weight $3$.

At the orthogonal intersection point between $A$ and $B$, there is an orbifold point of weight $18$ and local group $\bbZ / 3 \times \bbZ / 6$ generated by reflections through $A$ and $B$. Similarly, at the $\mathrm{A}_1$ singularity of $A$ one has orbifold weight $12$ and local group $\bbZ / 2 \times \bbZ / 3$, where the order $3$ generator is a reflection but the order $2$ generator is not and is associated with the singularity. The point $[1,e^{2 \pi i/3}, e^{4 \pi i/3}]$ in $\bbP^2$ projects to the $\mathrm{A}_2$ singularity of $\bbP(1,2,3)$, also marked with a $\times$, which has orbifold weight $3$ and group $\bbZ / 3$ that is not a reflection group.

Let $p[n]q$ denote the complex reflection group generated by elements of order $p$ and $q$ that satisfy a braid relation of order $n$. For example:
\begin{align*}
3[3]3 &= \langle a, b\ |\ a^3, b^3, a b a = b a b \rangle \\
3[4]6 &= \langle a, b\ |\ a^3, b^4, a b a b = b a b a \rangle
\end{align*}
The group $3[3]3$ is the exceptional Shephard--Todd group $G_4$ of order $24$, and this is the local orbifold group for the cusp singularity on $B$. The other point of intersection between $A$ and $B$ is a cusp of the orbifold (i.e., orbifold weight infinity), and one can check that the local group $3[4]6$ admits a central exact sequence
\[
1 \to \bbZ \to 3[4]6 \to \Del(2,3,6) \to 1,
\]
where $\Del(2,3,6)$ is the $(2,3,6)$ Euclidean triangle group. In other words, $\Del(2,3,6)$ is the index $2$ orientation-preserving subgroup of the group of Euclidean isometries generated by reflections in a triangle in the plane with angles $(\pi / 2, \pi / 3, \pi / 6)$. We note that one can consider the group $\bbZ / 3 \times \bbZ / 6$ associated with the orthogonal crossing of $A$ and $B$ as the complex reflection group $3[2]6$.

\medskip

The first author learned the following result from conversations with Domingo Toledo.

\begin{prop}\label{prop:DMPresent1}
Let $X$ be the Deligne--Mostow orbifold determined by the weights $(5,4,1,1,1)$ and $\Gam < \PU(2,1)$ be the associated lattice. Then $\Gam$ has presentation
\begin{equation}\label{eq:Pres1}
\Gam = \langle b, u, v\ |\ bv = vb, bub=ubu, uvuv=vuvu, (buv)^3, b^3, u^3, v^6\rangle.
\end{equation}
\end{prop}

\begin{pf}
One can check directly that the complement $\bbP(1,2,3) \ssm \{A, B\}$ has presentation with generators $b,u,v$ and relations
\[
\{bv = vb, bub=ubu, uvuv=vuvu, (buv)^3\},
\]
where $v$ is a loop around $A$ and $b, u$ denote loops around $B$ on either side of the cusp singularity (below and above in the configuration in Figure \ref{fig:DM}, respectively). The first three relations are braid relations arising from the two points of intersection between $A$ and $B$ and the cusp singularity of $B$. The fourth relation arises from the $\mathrm{A}_2$ singularity. To obtain a presentation for $\Gam$, one adds relations of the form $\mu^j$ where $\mu$ is a loop around the appropriate curve and $j$ is the orbifold weight of that curve. This is precisely the presentation in Equation \eqref{eq:Pres1}.
\end{pf}

\begin{rem}\label{rem:OrbifoldBraid}
The braid relations for the local groups $\bbZ / 3 \times \bbZ / 6$, $3[3]3$, and $3[4]6$ for the orbifold structure on $\bbP(1,2,3)$ are precisely the three braid relations (of orders two, three, and four, respectively) in the presentation for $\Gam$ from Proposition \ref{prop:DMPresent1}.
\end{rem}

It is known that $\Gam$ is arithmetic. One can directly compute its orbifold Euler characteristic, which is $1/72$, from Figure \ref{fig:DM}. It follows for example from \cite{StoverVol} that $\Gam$ is one of the two minimal covolume nonuniform arithmetic lattices in $\PU(2,1)$. In other words, it is a Picard modular group $\PU(h, \calO_k)$, where $k := \bbQ(e^{\pi i / 3})$, $\calO_k := \bbZ[e^{\pi i / 3}]$ is the ring of Eisenstein integers, and $h$ is a certain hermitian form on $k^3$.

For our purposes, it will be convenient to instead fix the hermitian form
\[
h_0 := \begin{pmatrix} 0 & 0 & 1 \\ 0 & 1 & 0 \\ 1 & 0 & 0 \end{pmatrix}.
\]
With respect to this basis, we will in fact see that $\Gam$ is the image in $\PU(2,1)$ of the subgroup of $\U(2,1)$ consisting of those matrices of the form
\begin{equation}\label{eq:IntegralForm}
\begin{pmatrix} a_{1,1} & a_{1,2} & a_{1,3} / \sqrt{-3} \\ \sqrt{-3}\, a_{2,1} & a_{2,2} & a_{2,3} \\ \sqrt{-3}\, a_{3,1} & \sqrt{-3}\, a_{3,2} & a_{3,3} \end{pmatrix}
\end{equation}
for $a_{j,k} \in \calO_k$. Generators are relations for this lattice were given by Zhao \cite{Zhao}. We will use round brackets to denote matrices in $\U(2,1)$ and square brackets to denote their image in $\PU(2,1)$. We have:

\begin{prop}\label{prop:DMPresent2}
With respect to the above hermitian form $h_0$, $\Gam$ is the subgroup of $\PU(2,1)$ generated by the two elements
\begin{align*}
x &:= \begin{bmatrix} 1 & 0 & -1/\sqrt{-3}\\ 0 & e^{5 \pi i/3} & 0 \\ \sqrt{-3} & 0 & 0 \end{bmatrix} \\
y &:= \begin{bmatrix} e^{5 \pi i/3} & e^{5 \pi i/3} & -1/\sqrt{-3} \\ \sqrt{-3} & e^{2 \pi i/3} & 0 \\ \sqrt{-3} & 0 & 0 \end{bmatrix}
\end{align*}
subject to the relations
\[
x^3 = y^3 = \left(y  x^{-1} y\right)^{12} = \left[\left(y x^{-1} y\right)^2\, ,\, \left(y x^{-1}\right)^2\right] = \Id,
\]
where $[\,\, ,\, ]$ denotes the commutator.
\end{prop}

\begin{pf}
This is the presentation for the subgroup of $\PU(2,1)$ consisting of those matrices satisfying Equation \eqref{eq:IntegralForm} given in \cite[Thm.\ 6.2]{Zhao}. One then checks either by hand or with Magma \cite{Magma} that the maps
\begin{align*}
x &\mapsto b \\
y &\mapsto b u v
\end{align*}
determine an isomorphism of this lattice with $\Gam$ that has inverse:
\begin{align*}
b &\mapsto x \\
u &\mapsto y x y^{-1} \\
v &\mapsto y x^{-1} y^{-1} x^{-1} y
\end{align*}
This implies the proposition.
\end{pf}

Finally, we describe the cusp subgroup of $\Gam$ in terms of our generators. According to \cite[Prop.\ 3.3]{Zhao}, the cusp subgroup has generators $r := (y x^{-1} y)^2$ and $s := y x^{-1}$ subject to the relations $r^6$, $(r^{-1} s)^3$, $[r, s^2]$. This is isomorphic to the group $3[4]6$ generated by $u$ and $v$ by the identifications:
\begin{align*}
u &\mapsto r s^{-1} \\
v &\mapsto r^{-1} \\
r &\mapsto v^{-1} \\
s &\mapsto u^{-1} v^{-1}
\end{align*}
We leave it to the reader to find matrix representatives for these generators using Proposition \ref{prop:DMPresent2}. We then have the following, which allows one to quickly verify Proposition \ref{prop:HirzebruchLattice}.

\begin{prop}\label{prop:HirzebruchinDM}
The lattice $\Gam_1 < \PU(2,1)$ defined by Hirzebruch's ball quotient $Z_1$ is a normal subgroup of index $72$ in $\Gam$.
\end{prop}

\begin{pf}
Since $Z_1$ has Euler characteristic $1$ and the Deligne--Mostow orbifold associated with $\Gam$ has orbifold Euler characteristic $1/72$, if $\Gam_1 < \Gam$, then it must be index $72$. In \cite{DiCerboStoverVol}, L.\ Di Cerbo and the first author showed that $\Gam$ is commensurable with $\Gam_1$, i.e., one can conjugate $\Gam_1$ in $\PU(2,1)$ so that its intersection with $\Gam$ is finite index in each.

To prove that $\Gam_1$ is in fact conjugate into $\Gam$, we will use Magma \cite{Magma}. One can show that there are two distinct homomorphisms
\[
\rho_1, \rho_2 : \Gam \to F,
\]
where $F$ is the finite group of order $72$ with identifier $\langle 72, 25\rangle$ in the Magma database of small finite groups. One can check that $F$ fits into a central exact sequence
\[
1 \to \bbZ / 6 \to F \to A_4 \to 1,
\]
where $A_4$ is the alternating group on four elements. We will show that $\Gam_1 = \ker(\rho_j)$ for one of these homomorphisms.

One such kernel $H$ has generators
\begin{align*}
h_1 &= x^{-1} y x y x y^{-1} x^{-1} y^{-1} \\
h_2 &= x^{-1} y^{-1} x^{-1} y x y x y^{-1} \\
h_3 &= x y^{-1} x y x y x^{-1} y^{-1} x \\
h_4 &= x y^{-1} x^{-1} y^{-1} x y x y x
\end{align*}
and relations:
\begin{align*}
h_3 h_2^{-1} h_1 h_4 h_2 h_1^{-1} h_3^{-1} h_4^{-1} &= \\
h_2 h_1^{-1} h_4^{-1} h_2^{-1} h_1 h_3 h_4 h_3^{-1} &= \\
h_2^{-1} h_1 h_3 h_1^{-1} h_3^{-1} h_4 h_2 h_4^{-1} &= \\
h_2 h_4 h_3^{-1} h_1^{-1} h_3 h_1 h_2^{-1} h_4^{-1} &= \\
h_1^{-1} h_3^{-1} h_2^{-1} h_4^{-1} h_3 h_4 h_1 h_2 &= \\ 
h_3^{-1} h_2^{-1} h_1 h_3 h_4^{-1} h_1^{-1} h_4 h_2 &= \\
h_1^{-1} h_3 h_1 h_3^{-1} h_4 h_2^{-1} h_4^{-1} h_2 &= \\
h_2^{-1} h_3 h_1 h_4^{-1} h_1^{-1} h_2 h_4 h_3^{-1} &= \\
h_3^{-1} h_4 h_2 h_4^{-1} h_3 h_4 h_3^{-1} h_2^{-1} h_3 h_4^{-1} &= \textrm{Id}
\end{align*}
In \cite[Prop.\ 6.3]{Zhao}, one finds a complete list of representatives for the conjugacy classes of finite order elements in $\Gam$. Checking in Magma that none of these elements is in $H$ implies that $H$ is torsion-free. Also, one checks that the image in $\Gam / H$ of the cusp subgroup $\langle r, s\rangle$ of $\Gam$ has index four. This implies that the associated ball quotient $\bbB^2 / H$ has four cusps.

We now check that the four cusps give a smooth toroidal compactification. Since $H$ is a normal subgroup of $\Gam$, it suffices to check that any one cusp gives a smooth toroidal compactification, as all four cusps are equivalent under the action of $F$ on the associated ball quotient. Indeed, the quotient is the Deligne--Mostow orbifold described above, which has one cusp. To see that one obtains a smooth toroidal compactification, one uses Magma. Recall that the cusp subgroup of $\Gam$ is generated by elements $r$ and $s$, defined above. Creating a homomorphism from $\langle r, s \rangle \cong 3[4]6$ to $\Gam/H$ and presenting the kernel, one sees that the kernel is isomorphic to the group
\[
\langle c_1, c_2, c_3\ |\ [c_1, c_2] c_3^{-1}, [c_1, c_3], [c_2, c_3] \rangle.
\]
This is precisely the cusp subgroup associated with a smooth toroidal compactification by an elliptic curve of self-intersection $-1$ (e.g., see \cite{DiCerboStoverVol}).

Then, $Z_1 := \bbB^2 / \Gam_1$ is the unique smooth ball quotient of Euler number one that admits a smooth toroidal compactification with four cusps \cite[Thm.\ 1.1]{DiCerboStoverVol}. We deduce that $H = \Gam_1$, completing the proof.
\end{pf}

\bibliography{HirzebruchCovers}

\end{document}